\documentclass[a4paper,dvips,12pt,oneside]{article}
\usepackage[makeroom]{cancel}
\usepackage{graphicx}
\usepackage{epstopdf}
\usepackage{subfig}
\usepackage{color}
\usepackage{float}
\usepackage{leftidx}
\usepackage{bigints}
\usepackage[affil-it]{authblk}
\usepackage{amssymb,amsmath,dsfont,url,epsfig}
\usepackage[left=1in, right=1in, top=1.1in, bottom=1.1in, includefoot, headheight=13.6pt]{geometry}
\usepackage[T1]{fontenc}
\usepackage{titlesec, blindtext, color}
\definecolor{gray75}{gray}{0.75}

\newcommand{\sln}{\linespread{1}}
\newcommand*{\email}[1]{\href{mailto:#1}{\nolinkurl{#1}} } 
\titleformat{\chapter}[block]{\LARGE\bfseries\sln}{Chapter \thechapter}{11pt}{\newline\huge\bfseries}
\newtheorem{theorem}{Theorem}[section]
\newtheorem{definition}{Definition}[section]
\newtheorem{example}{Example}[section]

\newtheorem{proposition}{Proposition}[section]

\begin{document}
\title{On Minkowskian Product Einstein-Finsler spaces}

\author[1]{Arti Sahu Gangopadhyay
  \thanks{E-mail: \texttt{arti.sahu@bhu.ac.in}}}
\affil[1,4]{DST-CIMS, Banaras Hindu University, Varanasi-221005, India}
\author[2]{Ranadip Gangopadhyay
  \thanks{E-mail: \texttt{gangulyranadip@gmail.com}}}
\affil[2]{Indian Institute of Science Education and Research (IISER) Mohali, Knowledge City, Sector 81, S.A.S. Nagar-140306, Punjab, India}
\author[3]{Ghanashyam Kr. Prajapati
  \thanks{E-mail: \texttt{gspbhu@gmail.com}}}
\affil[3]{Loknayak Jai Prakash Institute of Technology, Chhapra-841302, India}

\author[4]{Bankteshwar Tiwari
  \thanks{E-mail: \texttt{btiwari@bhu.ac.in}}}
\date{}
\maketitle

\begin{abstract}
    In this paper we study the  Minkowskian product Finsler manifolds. More precisely, we prove that if the Minkowskian product Finsler manifold is  Einstein then either the product manifold is Ricci flat or both the qoutient manifolds are Einstein with same scalar functions.\\
    

\end{abstract}
\section{Introduction}
The study of product manifolds in Riemannian and Finsler geometry has attracted the attention of differential geometers as well as physicists since long time. First reason of attraction is that it generates new spaces from the old and second is that it has lot of applications, specially in physics.  It is well known that there exists a canonical Riemannian metric on the product of two Riemannian manifolds, but as pointed out by Chern et al \cite{SSZ}, for a pair of Finsler manifolds $(\overline{M}, \overline{F})$ and $(\widetilde{M}, \widetilde{F})$, there is no canonical way to define Finsler metrics on the product manifold $ M = \overline{M} \times \widetilde{M}$. Hence, there are several methods to construct new Finsler metrics on $ M$, such as warped product, twisted product, Minkowskian product etc., see \cite{ CSZ, LK, TGS, TGS1}. 

Let $(\overline{M}, \overline{F})$ and $(\widetilde{M}, \widetilde{F})$ be two Finsler manifolds with Finsler metrics $\overline{F}, \widetilde{F}$, respectively, and $f: \overline{M} \times \widetilde{M} \to \mathbb{R}^+$ be a smooth function. In \cite{LK}, on the product manifold $M$ the authors consider the metric $F(y_1,y_2)=\sqrt{\overline{F}^2(y_1)+f^2(x_1,x_2)\widetilde{F}^2(y_2)}$ for all $(x_1,x_2)\in \overline{M}\times \widetilde{M}$ and $(y_1,y_2)\in {\overline{M}'\times \widetilde{M}'}$, where $\overline{M}'$ and $\widetilde{M}'$ are the slit tangent bundles $\overline{M}'= T\overline{M}\setminus \{0\}$, $\widetilde{M}'= T\widetilde{M}\setminus \{0\}$ respectively; and called it as twisted warped product metric. In \cite{CSZ}, Chen et al. defined the warped product Finsler metric $F=\overline{\alpha}\phi(s,r)$ on the product manifold $M=I\times \overline{M}$, where $(\overline{M}, \overline{\alpha})$ is a Riemannian manifold and $I$ is a interval in $\mathbb{R}$. They have also shown that this type of warped product of metrics contains the class of spherically symmetric Finsler metrics and contained in the class of general $(\alpha,\beta)$ Finsler metric. There is an another method to construct a Finsler metric of the form 
\begin{equation}\label{3.1}
 F(x_1,x_2,y_1,y_2):=\sqrt{f(\overline{\alpha}^2(x_1,y_1), \widetilde{\alpha}^2(x_2,y_2))}   
\end{equation}
on the product manifold $(\overline{M}\times \widetilde{M})$, where $(\overline{M},\overline{\alpha})$ and $(\widetilde{M},\widetilde{\alpha})$ are Riemannian manifolds, and $f(s, t): [0,\infty) \times [0,\infty) \to [0,\infty)$ be a $C^{\infty}$ function satisfies the following conditions:
\begin{enumerate}
\item  $f(\lambda s, \lambda t) = \lambda f(s,t), \forall \lambda >0$, and $f(s,t)>0, \forall(s,t) \ne (0,0)$.
\item $f_s>0, ~~ f_t>0, ~~ f_s+2sf_{ss}>0, ~~ f_t+2tf_{tt}>0$, and 
\item $f_sf_t-2ff_{st}>0$
\end{enumerate}
This Finsler metric $F$ defined in equation \eqref{3.1} is a non-Riemannian Berwald metric and affinely equivalent to the Riemannian metric $\sqrt{\overline{\alpha}^2+\widetilde{\alpha}^2}$, \cite{SSZ}. Here, $s=\overline{\alpha}^2,~ t= \widetilde{\alpha}^2$ and $f(s,t)=s+t$.\\

Minkowskian product of Finsler manifolds,  introduced by Okada \cite{OKD}, is a  generalisation of above mentioned product metrics.  He has proved that the Berwald connection of any Minkowskian product of Finsler spaces coincides with that of Euclidean product of these Finsler spaces. Further he has studied the geodesics of the Minkowskian product spaces. Later He, et al \cite{HLBN}, have derived Cartan connection and Berwald connection of Minkowskian product space in terms of the original Finsler spaces. They have also derived necessary and sufficients conditions for Minkowskian product space to be Berwald (resp. weakly Berwald, Landsberg, weakly Landsberg) spaces. Li et al \cite{LHZZ}, characterised dually flat (resp. projectively flat) Minkowskian product Finsler metrics and proved that Minkowskian product Finsler metric is dually flat and projectively flat Finsler metrics if and only if it is a Minkowskian metric. In this paper, at first we study the Riemann curvature and Ricci curvature of Minkowskian product of Finsler manifolds and obtain the following results:\\
\begin{proposition}\label{th3.1}
The Minkowskian product Finsler manifold $(M,F):=(\overline{M}\times \widetilde{M}, \overline{F}\times \widetilde{F})$ of the Finsler manifolds $(\overline{M},\overline{F})$ and $(\widetilde{M},\widetilde{F})$ has vanishing Riemann curvature if and only if $(\overline{M},\overline{F})$ and $(\widetilde{M},\widetilde{F})$ has vanishing Riemann curvature.    
\end{proposition}
\begin{theorem}\label{th3.2}
    The Minkowskian product Finsler manifold $(M,F):=(\overline{M}\times \widetilde{M}, \overline{F}\times \widetilde{F})$ of the Finsler manifolds $(\overline{M},\overline{F})$ and $(\widetilde{M},\widetilde{F})$  is Ricci-flat if and only if $(\overline{M},\overline{F})$ and $(\widetilde{M},\widetilde{F})$ are Ricci flat.   
\end{theorem}  
In the following theorem we study the Minkowskian product Einstein-Finsler metrics:
\begin{theorem}\label{th3.3}
The Minkowskian product Finsler manifold $(M,F):=(\overline{M}\times \widetilde{M}, \overline{F}\times \widetilde{F})$ of the Finsler manifolds $(\overline{M},\overline{F})$ and $(\widetilde{M},\widetilde{F})$  is Einstein, then, either
 $(M,F)$ is Ricci-flat or both $(\overline{M},\overline{F})$ and $(\widetilde{M},\widetilde{F})$ are Einstein with same scalar function. Moreover, in this case the product function $f$ is linear in $H$ and $K$, i.e., $f(H,K)=aH+bK$, where $a$ and $b$ are positive real constants.
\end{theorem}
\section{Preliminaries}
Let $ M $ be an $n$-dimensional smooth manifold. $T_{x}M$ denotes the tangent space of $M$
 at $x$. The tangent bundle $TM$ of $ M $ is the disjoint union of tangent spaces $T_xM$, i.e. $ TM:= \sqcup _{x \in M}T_xM $. We denote the elements of $TM$ by $(x,y)$, where $y\in T_{x}M $ and $TM_0:=TM \setminus\left\lbrace 0\right\rbrace $, the slit tangent bundle of $M$.\\
\begin{definition}
\textnormal{(\cite{SSZ}) A Finsler metric on $M$ is a function $F:TM \to [0,\infty )$ satisfying the following conditions:
 \begin{enumerate}
 \item $F$ is smooth on $TM\setminus \{0\}$,
  \item $F$ is positively 1-homogeneous on the fibers of the tangent bundle $TM$,
   \item The Hessian of $\displaystyle{\frac{F^2}{2}}$ with element $\displaystyle{g_{ij}=\frac{1}{2}\frac{\partial ^2F^2}{\partial y^i \partial y^j}}$ is nondegenerate on $TM_0$, where $y=y^i\frac{\partial}{\partial x^i}$.
  \end{enumerate}
 The pair $(M,F)$ is called a Finsler space and $g_{ij}$ is called the fundamental tensor.}
 \end{definition}
 Let $(\overline{M}, \overline{F})$ and $(\widetilde{M}, \widetilde{F})$ be two Finsler manifolds with dimension $m$ and $n$, respectively, then $M=\overline{M} \times \widetilde{M} $ is a product manifold with dimensions $m+n$. \\
 Let $(x^1,...,x^m)$ and $(x^{m+1},...,x^{m+n})$ be the local coordinates of $\overline{M}$ and $\widetilde{M}$, respectively, then the local coordinates on $M$ are $(x^1,...,x^{m+n})$. Let $(x^1,...,x^m,y^1,...,y^m)$ and $(x^{m+1},...,x^{m+n},y^{m+1},...,y^{m+n})$ be the induced local coordinates on the tangent bundle $T\overline{M}$ and $T\widetilde{M}$, respectively, then the induced local coordinates on the tangent bundle $TM$ are $(x^1,...,x^m,x^{m+1},...,x^{m+n},y^1,...,y^m,y^{m+1},...,y^{m+n})$. Denote $\overline{M}' = T\overline{M} \setminus\{0\}$, $\widetilde{M}'= T\widetilde{M} \setminus \{0\}$, $M'= \overline{M}' \times \widetilde{M}' \subset T(\overline{M} \times \widetilde{M}) \setminus \{0\}$. \\
  The following, lowercase Latin indices such as $i,j,k$, etc.,will run from $1$ to $m+n$, lowercase Latin indices such as $a,b,c$, etc.,will run from $1$ to $m$, whereas lowercase Greek indices such as $\alpha,\beta,\gamma$, etc.,will run from $m+1$ to $m+n$, and the Einstein summation convention is assumed throughout this paper.\\
  Let $f:[0,\infty ) \times [0,\infty ) \to [0,\infty )$ be a continuous function such that
  \begin{itemize}
      \item[(i)]  $f(s,t)=0$ if and only if $(s,t) = (0,0)$; 
      \item[(ii)] $f(\lambda s,\lambda t)=\lambda f(s,t)$ for any $\lambda \in [0,\infty )$;
      \item[(iii)] $f$ is smooth on $ (0,\infty ) \times (0,\infty )$;
     \item[(iv)] $\frac{\partial f}{\partial s} \neq 0$, $\frac{\partial f}{\partial t} \neq 0$ for any $ (s,t)\in (0,\infty ) \times (0,\infty )$;
      \item[(v)] $\frac{\partial f}{\partial s} \frac{\partial f}{\partial t}- 2f\frac{\partial ^2 f}{\partial s \partial t } \neq 0$ for any $  (s,t)\in (0,\infty ) \times (0,\infty )$.
 \end{itemize}

\begin{definition} \textnormal{ (\cite{OKD})}
\textnormal{Let $(\overline{M}, \overline{F})$ and $(\widetilde{M}, \widetilde{F})$ be two Finsler manifolds and $f$ be continuous function satisfying the conditions $\mathbf{(i)-(v)}$. Denote $K= \overline{F}^2$, $H=\widetilde{F}^2$, the Minkowskian product Finsler manifold $(M, F)$ of $(\overline{M},  \overline{F})$ and $(\widetilde{M}, \widetilde{F})$ with respect to the product function $f$, is the product manifold $M=\overline{M} \times \widetilde{M}$ endowed with the Finsler metric $F: M' \to \mathbf{R}^+$ given by
\begin{equation}
   F (x, y) = \sqrt{f (K (x^a, y^a),H (x^{\alpha}, y^{\alpha}))}, 
\end{equation}
where $(x, y) \in M'$, $(x^a, y^{a}) \in \overline{M}'$, $(x^{\alpha}, y^{\alpha}) \in \widetilde{M}'$ with $x =(x^a, x^{\alpha})$, $y=(y^a, y^{\alpha})$}.    
\end{definition}  
In $\cite{OKD}$, Okada showed that $F$  defined  above is a Finsler metric on $M$.
\begin{example}  \textnormal{ Let $f_1:[0,\infty ) \times [0,\infty ) \to [0,\infty )$ such that $f_1(s,t)=as+bt$, where $a,b$ are positive real numbers, then $f_1$  satisfies all the above conditions (i)-(v). Let $(\overline{M}, \overline{F})$ and $(\widetilde{M}, \widetilde{F})$ be two Finsler manifolds, then the Minkowskian product Finsler manifold $(M, F)$ of $(\overline{M},  \overline{F})$ and $(\widetilde{M}, \widetilde{F})$ with respect to the product function $f_1(s,t)=as+bt$ is given by $F (x, y) = \sqrt{a \overline{F}^2 (x^a, y^a)+b\widetilde{F}^2(x^{\alpha}, y^{\alpha}))}, $ which is clearly a Finsler manifold.
} 
\end{example}
\begin{example}  \textnormal{ Let   $f_2:[0,\infty ) \times (0,\infty ) \to [0,\infty )$ such that $f_2(s,t)=\frac{s^2}{t}$, then  $f_2$ also satisfies all the above conditions (i)-(v). Let $\overline{F}$ be a Randers metric given by $\overline{F}=\alpha + \beta$ where, $\alpha$ be a Riemannian metric and $\beta$ be a 1-form with $\|\beta\|_{\alpha}<1$ and $\widetilde{F}$ be a Riemannian metric given by $\widetilde{F}=\alpha$. Then for $ f_2= \frac{s^2}{t}$, we get Minkowskian product Finsler metric $F=\frac{(\alpha +\beta)^2}{\alpha}$ of Randers metric $\overline{F}$ and Riemannian metric $\widetilde{F}$, which is well known Z. Shen's square metric. } 
\end{example}
\begin{proposition}\textnormal{(\cite{HLBN})}
Let $(M,F)$ be a Minkowskian product of the Finsler manifolds $(\overline{M},\overline{F})$ and $(\widetilde{M},\widetilde{F})$. Then the fundamental tensor matrix of $F$ is given by 
\begin{equation}
 g_{ij}= \left( \frac{\partial^2 F^2}{\partial y^i \partial y^j}\right)= \begin{pmatrix} G_{ab} & G_{a\beta} \\ G_{\alpha b} & G_{\alpha \beta}, \end{pmatrix} 
\end{equation}
where
\begin{equation}\label{eq3.3}
\begin{split}
G_{ab}= f_KK_{ab}+f_{KK}K_aK_b, \qquad G_{a\beta}= f_{KH}K_aH_{\beta},\\ G_{\alpha b}= f_{KH}H_{\alpha}K_b, \qquad  G_{\alpha \beta}=f_HH_{\alpha\beta}+f_{HH}H_{\alpha}H_{\beta}.   
\end{split}    
\end{equation}
and
\begin{equation*}
 K_i=\frac{\partial K}{\partial y^i}, ~ K_{;i}=\frac{\partial K}{\partial x^i}, ~ K_{i;j}=\frac{ \partial ^2K}{ \partial y^i\partial x^j},   
\end{equation*}
\begin{equation*}
 H_i=\frac{\partial H}{\partial y^i}, ~ H_{;i}=\frac{\partial H}{\partial x^i}, ~ H_{i;j}=\frac{ \partial ^2H}{ \partial y^i\partial x^j},   
\end{equation*}
\end{proposition}
\begin{proposition}\textnormal{(\cite{HLBN})}
Let $(M,F)$ be a Minkowskian product of the Finsler manifolds $(\overline{M},\overline{F})$ and $(\widetilde{M},\widetilde{F})$. Then the inverse matrix of the fundamental tensor matrix of $F$ is given by 
\begin{equation}
 g^{ij}= \begin{pmatrix} G^{ba} & G_{b\alpha} \\ G_{\beta a} & G_{\beta \alpha}, \end{pmatrix} 
\end{equation}
where
\begin{equation}
\begin{split}
G^{ba}=\frac{1}{f_K}\left( K^{ba}- \frac{f_Hf_{KK}}{\Delta}y^by^a\right), \qquad G_{b\alpha}= -\frac{1}{\Delta}f_{KH}y^by^{\alpha},\\ G_{\beta a}= -\frac{1}{\Delta}f_{KH}y^{\beta}y^a, \qquad  G_{\alpha \beta}= \frac{1}{f_H}\left( H^{\beta\alpha}- \frac{f_Kf_{HH}}{\Delta}y^{\beta}y^{\alpha}\right).   
\end{split}    
\end{equation}
and $\Delta=f_Kf_H-2ff_{KH}$
\end{proposition}
\begin{definition} \textnormal{(\cite{ShiBanktesh})
A smooth curve in a
Finsler space is a geodesic if it has constant speed and is locally length
minimizing. Thus a geodesic in a Finsler space $(M,F)$ is a curve
$\gamma:I=[a,b]\rightarrow M$ with $F(\gamma (t), \dot{\gamma }(t))=$ constant and
for any $t_0\in I$, there is a small number $\epsilon > 0$ such that
$\gamma$ is length minimizing on
$[t_0-\epsilon,t_0+\epsilon]\cap I$.}
\end{definition}
It can be shown that a smooth curve $\gamma$ in a Finsler manifold $(M,F)$ is a geodesic if and only if $\gamma(t)$
satisfies the following:
\begin{equation}\label{2.2a}
\frac{d^2\gamma^i(t)}{dt^2}+{G}^i\left(\gamma,\frac{d \gamma}{dt}\right)=0,
\end{equation}
where $G^i=G^i(x,y)$ are local functions on $TM$ defined by
\begin{equation}\label{eq1.2}
{G}^i=\frac{1}{4}{g}^{i{\ell}}\left\{\left[{F}^{2}\right]_{x^ky^{\ell}}y^k-
\left[{F}^{2}\right]_{x^{\ell}}\right\}.
\end{equation}
The coefficients $G^i$ are called the spray coefficients.
\begin{theorem} \textnormal{(\cite{OKD})}
  The Spray coefficients of the Minkowskian product $(M,F)$ of the Finsler manifolds $(\overline{M}, \overline{F})$ and $(\widetilde{M}, \widetilde{F})$ are 
  \begin{equation}
     G^{i}= \begin{cases}\overline{G}^a  \quad : i = a \\  \widetilde{G}^{\alpha}   \quad :  i = \alpha\end{cases}
  \end{equation}
  where $\overline{G}^a$ and $\widetilde{G}^{\alpha}$ are spray coefficients of the manifolds $(\overline{M}, \overline{F})$ and $(\widetilde{M}, \widetilde{F})$ respectively.
\end{theorem}
 \section{Main results}
 In this section we study some non-Riemannian curvature properties of Minkowskian product of two Finsler metrics and prove the theorems mentioned in the Introduction.
 \subsection{Riemann curvature}
 \begin{definition}\label{def2}
  \textnormal{The Riemann curvature
$R={{R}_y:T_xM \rightarrow T_xM}$, for a Finsler space ${F^n,}$ is defined by
${R}_{y}(u)={R}^{i}_{k}(x,y)u^{k} \frac{\partial}{\partial x^i}$,
\ $u=u^k\frac{\partial}{\partial x^k}$,\ where ${R}^{i}_{k}={R}^{i}_{k}(x,y)$ denote the coefficients of the Riemann curvature of $F$ and given by
\begin{equation}\label{eqn3.1}
{R}^{i}_{k}=2\frac{\partial{G}^i}{\partial
x^k}-y^j\frac{\partial^2{G}^i}{\partial x^j\partial
y^k}+2{G}^j\frac{\partial^2 {G}^i}{\partial y^j\partial
y^k}-\frac{\partial{G}^i}{\partial y^j}\frac{\partial
{G}^j}{\partial y^k}.
\end{equation}}
\end{definition}
Let us denote the coefficients of Riemann curvature of $M$, $\overline{M}$, $\widetilde{M}$ by $R^{i}_{j}$, $\overline{R}^{a}_{b}$, $\widetilde{R}^{\alpha}_{\beta}$.\\

\textbf{Proof of Proposition \ref{th3.1}} Let us assume the coefficients of Riemann curvature tensors of $M$, $\overline{M}$ and $\widetilde{M}$ are denoted by $R^i_k$, $\overline{R}^a_b$ and $\widetilde{R}^{\alpha}_{\beta}$ respectively.\\

 Let  $i = a, k=b$,  therefore,
\begin{equation*}
   {R}^{i}_{k} =2\frac{\partial{G}^a}{\partial
x^b}-y^c\frac{\partial^2{G}^a}{\partial x^c\partial
y^b}+2{G}^c\frac{\partial^2 {G}^a}{\partial y^c\partial
y^b}-\frac{\partial{G}^a}{\partial y^c}\frac{\partial
{G}^c}{\partial y^b}
-y^{\beta}\frac{\partial^2{G}^a}{\partial x^{\beta}\partial
y^b}+2{G}^{\beta}\frac{\partial^2 {G}^a}{\partial y^{\beta}\partial
y^b}-\frac{\partial{G}^a}{\partial y^{\beta}}\frac{\partial
{G}^{\beta}}{\partial y^b}
\end{equation*}
Since $G^a$ is independent of $y^{\beta}$  implies , 
 $\frac{\partial{G}^a}{\partial
y^{\beta}}=0$,  therefore
\begin{equation*}
   {R}^{i}_{k} = \overline{R}^{a}_{b}
\end{equation*}
Similarly, if $i= \alpha , k = \beta$, then ${R}^{i}_{k} = \widetilde{R}^{\alpha}_{\beta}$.\\
Then
if $ i= a , k=\beta $, then 
\begin{equation}
 {R}^{i}_{k} = 2\frac{\partial{G}^a}{\partial
x^{\beta}}-y^c\frac{\partial^2{G}^a}{\partial x^c\partial
y^{\beta}}+2{G}^b\frac{\partial^2 {G}^a}{\partial y^b\partial
y^{\beta}}-\frac{\partial{G}^a}{\partial y^b}\frac{\partial
{G}^b}{\partial y^{\beta}}
-y^{\gamma}\frac{\partial^2{G}^a}{\partial x^{\gamma}\partial
y^{\beta}}+2{G}^{\gamma}\frac{\partial^2 {G}^a}{\partial y^{\gamma}\partial
y^{\beta}}-\frac{\partial{G}^a}{\partial y^{\gamma}}\frac{\partial
{G}^{\gamma}}{\partial y^{\beta}}   
\end{equation}
Since $G^{a}$ is independent of $y^{\beta}$ and $x^{\beta}$ implies , $\frac{\partial{G}^a}{\partial
x^{\beta}}=0$, $\frac{\partial{G}^a}{\partial
y^{\beta}}=0$, therefore
\begin{equation}
    {R}^{a}_{\beta}=0
\end{equation}
Then
\begin{equation}\label{eq3.8}
R^i_k= \begin{cases}\overline{R}^a_b \quad : i=a,k = b \\  \widetilde{R}^{\alpha}_{\beta} \quad : i= \alpha,k=\beta \\ 0 \qquad :\textnormal{otherwise} \end{cases}
\end{equation}
Let us suppose $(M,F)$ has vanishing Riemann curvature. Then $R^i_k=0$. Hence, from \eqref{eq3.8} we get $\overline{R}^a_b=0$ and $\widetilde{R}^{\alpha}_{\beta}=0$. Hence, Riemann curvature of $\overline{M}$ and $\widetilde{M}$ vanishes.\\
Conversely suppose Riemann curvature of $\overline{M}$ and $\widetilde{M}$ vanishes. Hence, $\overline{R}^a_b=0$ and $\widetilde{R}^{\alpha}_{\beta}=0$. Therefore, from \eqref{eq3.8} it is immediate that $R^i_k=0$ for all $i,k$. Hence, $(M,F)$ has vanishing Riemann curvature. 

 \subsection{Ricci curvature}

 \begin{definition}\label{def3} \textnormal{The trace of Riemann curvature $R^i_k$ is called the Ricci curvature, denoted by \textit{Ric(x,y)}  and defined by
 \begin{equation}
   Ric(x,y)=R^i_i(x,y).
 \end{equation}}\end{definition}
 The Ricci curvature $Ric$ is a scalar function on $TM_{0}$.
 \begin{definition}
     A Finsler metric is called an Einstein metric if $Ric = (n-1)cF^2$ for some scalar function $c = c(x)$ on $M$. 
 \end{definition}
 \begin{definition}
     The Ricci tensor of F is defined by 
     \begin{equation}
         Ric_{ij}:=\left[\frac{1}{2}Ric\right]_{y^iy^j}=\frac{1}{2}[R_s^s]_{y^iy^j}
     \end{equation}
 \end{definition}
\textbf{Proof of Theorem \ref{th3.2}}
Let us assume the Ricci curvature of $M$, $\overline{M}$ and $\widetilde{M}$ are denoted by $Ric$, $\overline{Ric}$ and $\widetilde{Ric}$ respectively. Then 
\begin{equation}\label{eq3.5}
Ric= \sum\limits_{k=1}^{m+n}{R}^{k}_{k }= \sum\limits_{k=1}^{m}{R}^k_k +\sum\limits_{k=n+1}^{m+n} {R}^{k}_{k}=\sum\limits_{a=1}^{m}\overline{R}^a_a +\sum\limits_{\alpha=m+1}^{m+n} \widetilde{R}^{\alpha}_{\alpha}=\overline{Ric}+\widetilde{Ric}
\end{equation}
Now we suppose $(M,F)$ is Ricci flat. Then $Ric=0$. Therefore, $ \overline{Ric}+\widetilde{Ric}=0$ i.e. $\overline{Ric} = - \widetilde{Ric}$. Since, $\overline{Ric}$ and $\widetilde{Ric}$ are independent of each other. Therefore, $\overline{Ric} = \widetilde{Ric}=0$. Hence, $(\overline{M}, \overline{F} )$ and $(\widetilde{M},\widetilde{F})$ are Ricci flat.
Therefore, if $\overline{Ric} = 0$ and  $\widetilde{Ric} = 0$ then clearly $Ric= 0$.
\par Conversely suppose $(\overline{M}, \overline{F} )$ and $(\widetilde{M},\widetilde{F})$ are Ricci flat. Then $\overline{Ric} =0$ and $\widetilde{Ric} = 0$. Therefore, $Ric=\overline{Ric}+\widetilde{Ric}=0$. Hence, $(M,F)$ is Ricci flat. 
\begin{definition}
 \textnormal{An $n$-dimensional Finsler metric $F$ is called an Einstein metric if its Ricci curvature $Ric$ is isotropic, i.e., $Ric_{ij} = (n - 1)\lambda g_{ij}$, where $\lambda = \lambda(x)$ is a scalar function and $Ric_{ij}=[Ric]_{y^iy^j}$.}     
\end{definition}
The Finsler metric $F$ is said to be of Ricci constant, if $\lambda$ = constant. In particular, $F$ is said to be Ricci-flat if $\lambda = 0$.\\

\textbf{Proof of Theorem \ref{th3.3}} Let $(M,F)$ be a Einstein metric. Then 
\begin{equation}\label{eq3.1}
Ric_{ij}=(n-1)\lambda(x)g_{ij}
\end{equation}
Now let $i\in \{1,...,m\}$ and $j\in \{m+1,...,m+n\}$, i.e., $i=a,b$ and $j=\alpha,\beta$. Then  differentiating \eqref{eq3.5} first with respect to $y^a$ and then with respect to $y^{\beta}$ we get
\begin{equation}
Ric_{a\beta}=   \overline{Ric}_{a\beta}+\widetilde{Ric}_{a\beta} 
\end{equation}
Since, $\overline{Ric}$ is independent of $y^{\beta}$, we have $\overline{Ric}_{a\beta}=0$. Similarly since $\widetilde{Ric}$ is independent of $y^{b}$, we have $\widetilde{Ric}_{a\beta}=0$   Hence, \eqref{eq3.1} yields 
\begin{equation}\label{eq3.2}
  (n-1)\lambda(x)g_{a\beta}=0 \qquad  (n-1)\lambda(x)g_{\alpha b}=0.
\end{equation}
From the first equation of \eqref{eq3.2} either $\lambda(x)=0$ or $g_{a\beta}=0$. Now consider the above cases:\\
\textbf{Case-1:} If $\lambda(x)=0$, then from \eqref{eq3.1} we obtain $Ric_{ij}=0$. \\
\textbf{Case-2:} If $\lambda(x)\ne 0$, then from \eqref{eq3.2}, we get $g_{a\beta}=0$ and $g_{\alpha b}=0$. Hence, from \eqref{eq3.3} we obtain
\begin{equation*}
    G_{a\beta}= f_{KH}K_aH_{\beta}=0,\\ G_{\alpha b}= f_{KH}H_{\alpha}K_b=0,
\end{equation*}
Now we consider the following cases:\\
\textbf{subcase-1}: $f_{KH}=0$. Then $f=aH+bK$.\\
Now differentiating \eqref{eq3.5} first with respect to $y^a$ and then with respect to  $y^b$  we get
\begin{equation}
Ric_{ab}= \overline{Ric}_{ab}
\end{equation}
Here we used the fact that $\widetilde{Ric}_{ab}=0$. Now from \eqref{eq3.1} we get 
\begin{equation*}
\overline{Ric}_{ab}=(n-1)\lambda(x)G_{ab}.    
\end{equation*}
Hence, $(\overline{M},\overline{F})$ is a Einstein metric.\\
By similar arguments $(\widetilde{M}, \widetilde{F})$ is also Einstein.\\
\textbf{subcase-2}: $f_{KH}\ne 0$. Then either $K_a=0$, i.e., $K$ is Riemannian, or, $H_{\alpha}=0$, i.e., $H$ is Riemannian.\\
If $K_a=0$, then from first, second and third equality of \eqref{eq3.3} $g_{ij}$ becomes a singular matrix, which is a contradiction.  Similarly $H_{\alpha}=0$ also gives a contradiction.

\subsection*{Statements \& Declarations}
The first author is supported by UGC Non-Net Fellowship by Banaras Hindu University second author is supported by
IISER Mohali Institute Post-doctoral fellowship.

 \end{document}